\documentclass[11pt,fleqn,a4paper]{article} 
\usepackage{amsfonts,amsmath}
\usepackage[english]{babel}
\usepackage[T1]{fontenc}
\usepackage{dsfont}
\usepackage[latin1]{inputenc}
\usepackage{a4}
\usepackage{mathbbol}
\usepackage{amsfonts,amssymb,amsmath,epsfig}
\usepackage{bbold}
\usepackage{stmaryrd,epsfig}
\usepackage{pifont,fancyheadings,graphicx,subfigure,ifthen,epsfig}
\usepackage{enumerate,rotating,hangcaption} 
\newcommand{\ep}{\epsilon}

\newcommand{\w}{\omega} 

\newcommand{\h}{\hat}

\newcommand{\Dx}{{\partial_x}} 
\newcommand{\Dy}{{\partial_y}} 
\newcommand{\B}{\bar} 
\newcommand{\Dxx}{{\partial^2_{xx}}} 
\newcommand{\Dyy}{{\partial^2_{yy}}} 
\newcommand{\Dxy}{{\partial^2_{xy}}} 
\newcommand{\dx}{{\partial_x}} 
\newcommand{\dy}{{\partial_y}} 
 
\newcommand{\3}{{ {\mathds 1}}_{x \geq 0}} 
\newcommand{\5}{{ {\mathds 1}}_{x \leq 0}} 
\newcommand{\4}{{ {\mathds 1}}_{y \geq 0}} 
 
\newcommand{\8}{{ {\mathds 1}}_{y \leq 0}}

\newcommand{\R}{\mathds R} 
\newcommand{\C}{\mathds C} 
\newtheorem{theoreme}{Theorem} 
\newtheorem{Def}{Definition}
\newtheorem{proposition}{Proposition} 
\newtheorem{lemma}{Lemma} 
\newtheorem{coroll}{Corollary}
\newenvironment{proof}{{\bf Proof.}}{$\diamondsuit$} 

\begin{document}

 
\title{Boundary Value Problem for an Oblique Paraxial Model of Light Propagation}

\author{Marie Doumic \thanks{INRIA Rocquencourt, projet BANG, Domaine de Voluceau, BP 105, F 78153 Rocquencourt, France; email: marie.doumic@inria.fr}}
 
\maketitle 

\begin{abstract}
We study the Schr\"odinger equation which comes from  the paraxial approximation of 
the Helmholtz equation in the case where the direction of propagation is tilted with 
respect to the boundary of the domain.  This model has been proposed in \cite{DGS}. Our primary interest here is in the boundary conditions successively in a half-plane, then in a quadrant of $\R^2.$  The half-plane problem has been used in \cite{DDGS} to build a numerical method, which has been introduced in the $HERA$ plateform of $CEA.$ 
\end{abstract}







                                \section{Introduction}\label{intro-2}

The simulation of a laser wave according to the paraxial approximation of the Maxwell equation has been intensively studied for a long time when the simulation domain is rectangular and the direction of propagation is parallel to one of the
principal axis of simulation domain(see for instance \cite{fcs} and references therein). 

We are concerned here with the case where the direction of
propagation is \emph{not} parallel to a principal axis of the simulation domain, and cannot even be considered as having a small incidence angle with it. It may be crucial for example if we want to simulate several beams with different directions,
possibly crossing each other in the same domain. This tilted frame model has been considered some years ago by physicists for dealing with beam crossing problems (see \cite{tik}), and is of particular interest in the framework of CEA's Laser Megajoule experiment (see \cite{sen, ria}). 

In \cite{DGS, DDGS}, a new model was proposed to deal with this case, and a numerical scheme was introduced and coupled with a time-dependent interaction model (this scheme was then used in the \emph{HERA} plateform of CEA). We focus here on the theoretical study of this new model, of ``advection-Schr\"odinger'' type.

This model is derived from the paraxial approximation, intensively used in optic models and in a lot of models related to laser-plasma interaction in Inertial Confinement Fusion (cf \cite{dorr, faf-2, foc-2, feu, para-3, rose} and their references). According to laws of optics, the laser electromagnetic field may be modeled by  the following frequency Maxwell equation (the Helmholtz problem):
\begin{equation}
\epsilon ^{2}\mathbf{\Delta }\psi +\psi +i2\nu_t \epsilon \psi =0, \qquad \overrightarrow{x}\;\in\;{\cal D}
\label{base}
\end{equation}
\ where $\epsilon ^{-1}$ is the wave number of the laser wave in a medium corresponding to the mean value of the refractive index and  $\nu_t =\nu + i \mu$
is a complex coefficient : its real part $\nu$ corresponds to a conveniently scaled absorption coefficient and its imaginary part $\mu$ to the variation of the refraction index.

The first basic assumption is that the wave length $2\pi
\epsilon $ is small compared to the size of the simulation domain (indeed, it is in the order of $1 \mu m$ and the simulation domain in in the order of some $mm$ for the Inertial Confinement plasmas). Assume also 
 that the light propagates according a fixed direction defined by the unit vector $\overrightarrow{k}.$
 Let us denote the longitudinal and the transverse coordinates by $z=\overrightarrow{x}\cdot\overrightarrow{k}$ and $X=\overrightarrow{x}-(\overrightarrow{x}\cdot\overrightarrow{k})
\overrightarrow{k}$ ; the gradient with respect to $ X$  is denoted $\nabla _{\bot } = \nabla - \overrightarrow{k}(\overrightarrow{k}.\nabla ) $. Now, if we replace $\psi $ by $u\exp (\frac{i\overrightarrow{k}\overrightarrow{x}}{\epsilon }),$ Equation (\ref
{base}) may be written as:
\begin{equation*}
\epsilon (2i\nu_t u+2i\overrightarrow{k}\cdot\overrightarrow{\mathbf{\nabla }}%
u)+\epsilon ^{2}\Delta _{\bot }u=-\epsilon ^{2}\frac{\partial ^{2}}{%
\partial z^{2}}u.
\end{equation*}
where $\Delta _{\bot }$ is the Laplacian operator in the transverse variable.
Assuming that $u$ is slowy varying with respect to the longitudinal
variable, we can neglect the right-hand side of the previous equation (we will check this \emph{a posteriori} ; see also \cite{D} for more details).
Therefore $u$ satisfies the classical paraxial  equation for wave propagation, which is a linear Schr\"{o}dinger equation:
\begin{equation}
i\overrightarrow{k}\cdot\overrightarrow{\mathbf{\nabla }}u+\frac{\epsilon }{2}%
\Delta _{\bot }u+i\nu_t u=0,\qquad \overrightarrow{x}\;\in\;{\cal D}. \label{bbb}
\end{equation}
If the laser beam enters the simulation domain with a very small incidence angle that is, say, if the vector $\overrightarrow{k}$  is normal to a boundary  of the simulation
domain, the boundary condition on the entrance boundary $\{x=0\}$  is straightforward: it is simply the value of the solution, which is a data. On the face $\{y=0\}$ of the domain parallel to the vector $\overrightarrow{%
k},$
to  deal correctly with the boundary
conditions 
 one has to introduce a fractional
derivative according to the $x$ variable - see \cite{Arnold, Arnold-2, DM, DM-2, joly} for a mathematical approach and \cite{had-2} for a physical approach ; see also Appendix \ref{racine-2}, or \cite{bcd} for numerical treatment.

But, when the incidence angle between $\overrightarrow{k}$ and the normal vector to the entrance boundary $\overrightarrow{e_x}$ is not small any more, what is of particular interest for the treatment of beam-crossing for instance, the classical approach is no more valid.

The aim of this article is to analyse Equation (\ref{bbb}) with an arbitrary angle between  $\overrightarrow{k}$ and  $\overrightarrow{e_x}$, first in a half-space, and then in a quadrant, and to address the issue to find proper boundary conditions on the edges of the domain.
The half-space problem has been studied and then used numerically in \cite{DDGS}; see also \cite{desr} for results of numerical simulations. 
To solve it, we first have to formulate the right entrance boundary condition and to show that the corresponding problem is well-posed.
For this purpose, we consider a 2D problem but most of the ideas developed here may be extended to the 3D case.
 For the statement of the entrance boundary condition,
 one assumes that a fixed plane wave $\psi^{in} = u^{in}\exp (\frac{i
\overrightarrow{k}\overrightarrow{x}}{\epsilon })$ is coming into the domain where $u^{in}$ is a given function of the variable $y$. Remark that for the Helmholtz problem, the boundary condition  is classical and may read as:
$$
(\epsilon \frac{\partial }{\partial \overrightarrow{n}}+i\overrightarrow{k%
}\cdot\overrightarrow{n}) (\psi -u^{in}e^{i\overrightarrow{k}\overrightarrow{x}%
\mathbf{/\epsilon }})=0,\qquad \overrightarrow{x}\;\in\;\partial{\cal D}. 
$$
where $\overrightarrow{n}$ is the outwards unit normal to the domain. By using the definition
$$\psi = u\exp (\frac{i \overrightarrow{k} .\overrightarrow{x}}{\epsilon })$$ in this expression,  the corresponding entrance boundary condition for Equation (\ref{bbb}) may be written  in a natural way as: 
\begin{equation}\label{Cbord}
(\epsilon \overrightarrow{n}\cdot \nabla _{\bot } +2i\overrightarrow{k}\overrightarrow{n}) (u-u^{in})=0,\qquad \overrightarrow{x}\;\in\;\partial{\cal D}.
\end{equation}
 In Section \ref{demi-2}, for the sake of completeness and clarity, we recall (see \cite{DDGS}, Section 2) and detail how to analyse Problem (\ref{bbb})(\ref{Cbord}) in the half-space $${\cal D}=\{(x,y) \quad\emph{s.t.}\quad x\geq 0\},$$ which is the simplest case. First we recall a classical energy estimate in the space $L^{2}(\R_+\times \R),$ which implies uniqueness of the solution for any $\nu_t $ satisfying \begin{equation}\label{hyp1}
\hbox {inf}_{\overrightarrow{x}} \nu(\overrightarrow{x})  > 0.
\end{equation}
Then, in the case where $\nu$ is a stritly positive constant, we recall in Theorem \ref{demi-theo-2} the exact analytical formulation of the solution, which is central in \cite{DDGS} to build a numerical scheme, based on a splitting technique and on Fast Fourier Transform. We also develop here some corollaries on the regularity of the solution.

We then have the tools to focus, in Section \ref{quart}, on Problem (\ref{bbb})(\ref{Cbord}) in the quadrant 
$${\cal D}=\{(x,y) \emph{ s.t.} \ x\geq 0,y\geq 0\},$$
which will make large uses of the solution on the half-space.

The most delicate point is to find a proper condition on the boundary $\{y=0\},$ which would be such that the restricted problem:
\begin{enumerate}
\item is well-posed,
\item admits for a solution the restriction to the quadrant of the solution in the half-space.
\end{enumerate}
Such a condition is called a \emph{transparent boundary condition}. If the solution of the quadrant is not the exact restriction of the solution in the half plane, but is close to it in some sense, the boundary condition is called an \emph{absorbing boundary condition}.

For reasons which will appear in the sequel, the proper choice for such a condition is given by Equation (\ref{CT}), which writes:
$$ 
\Dy U_{|y=0} - A_+ (\Dx) (U_{|y=0}) = 0 \;\; \forall x>0,$$
with $A_+$ a pseudo-differential operator  defined by
$$A_+ (\Dx)= {k_y \over k_x} \Dx -i {k_y\over \ep k_x^2}\bigl(1 + \sqrt{1+{2i\ep k_x\over k_y^2} \Dx +2i\ep\nu {k_x^2\over k_y^2}}\bigr).$$
For a definition of fractional derivatives, see Appendix \ref{racine-2} or \cite{Schwartz}, ch. VI.5. 
The derivation of this formula is similar to the derivation of the transparent boundary condition given for instance in \cite{Arnold, Arnold-2, DM, Baskakov} for the Schr\"odinger equation, using Fourier-Laplace transform: here, the operator $\overrightarrow{k}\cdot \overrightarrow{\nabla}$ plays the role of time, and if we suppose $k_y=0$ we would recover a ``time'' fractional derivative as in the previous references. 

Theorem \ref{CT-extension} then reformulates Problem (\ref{bbb})(\ref{Cbord})(\ref{CT}) for functions whose derivatives may not be defined at the boundary. The main result of this paper is then given by Theorem \ref{quart-theo1}: it states well-posedness and an explicit formulation for the solution to Problem (\ref{bbb})(\ref{Cbord})(\ref{CT}) and analyses its relation with the half-space solution, showing that if the ray enters in the quadrant (\emph{i.e} if $\overrightarrow{k}\cdot \overrightarrow{n} \leq 0$ ) Condition (\ref{CT}) is a \emph{transparent} one, and only an \emph{absorbing} one if it goes out of it (\emph{i.e} if $\overrightarrow{k}\cdot \overrightarrow{n} \leq 0$).

\section{Recall on the Half-Space Problem} \label{demi-2}

Here we address a problem where the simulation domain is the half-space 
$${\cal D}=\{\overrightarrow{x}\cdot \overrightarrow{e_x} \geq 0\}.$$
Set $ \overrightarrow{x}=(x,y)$ and define $\overrightarrow{k}=(k_x, k_y)$ to be a unitary vector of the plane which indicates the direction of the laser beam (we must of course have $k_x >0$). Assuming that (\ref{hyp1}) holds and that $\mu$ is a bounded function,  we rewrite Problem (\ref{bbb})(\ref{Cbord}) under the form 
\begin{equation}
i\overrightarrow{k}\cdot\overrightarrow{\mathbf{\nabla }}u+\frac{\epsilon }{2}%
\Delta _{\bot }u- \mu u+i\nu u =0, \label{bbbb}
\end{equation}
\begin{equation}\label{CE0}
(i\epsilon D-2 \overrightarrow{k}.\overrightarrow{n} ) (u-u^{in})|_{x=0}=0,  \quad  \end{equation}
where we have denoted $D=\overrightarrow{n}\cdot \nabla _{\bot }=k_{y}(k_{x}\partial _{y}-k_{y}\partial
_{x})$.

We first recall (See \cite{DDGS}, Proposition 1.) that if $u \in H^1(\R_+\times \R)$ is a solution to Problem (\ref{bbbb}) (\ref{CE0}) in the half-plane $\{x\geq 0\},$ we have the two following equivalent identities :

\begin{equation*}
\begin{array}{l} \displaystyle{
\iint\limits_{\mathcal{D}}2\nu |u |^{2}
+\int\limits_{\Gamma_0}\frac{|\overrightarrow{k}\cdot\overrightarrow{n}|}{2}|  
\frac{(i\epsilon D + 2 \overrightarrow{k}\cdot\overrightarrow{n} )u}{2|\overrightarrow{k}\cdot\overrightarrow{n}|}| ^2
 = }
\displaystyle{\int\limits_{\Gamma_0 }  |\overrightarrow{k}\cdot\overrightarrow{n}|\biggl(|u|^2 
+\frac{1}{2}  
|\frac{(i\epsilon D  - 2 \overrightarrow{k}\cdot\overrightarrow{n}) u^{in} }{2|\overrightarrow{k}\cdot\overrightarrow{n}|}| ^2
\biggr),}
\\ \\
\displaystyle{ 
 \iint\limits_{\mathcal{D}}2\nu |u |^{2} + \int\limits_{\Gamma_0}  |\overrightarrow{k}\cdot\overrightarrow{n}| |u|^2  = -{\cal I}m \bigl( \int\limits_{\Gamma_0} \B{u} (\ep D + 2i\overrightarrow{k}\cdot\overrightarrow{n}){u^{in}} \bigr).}
\end{array}
\end{equation*}
The first estimate can be interpreted as the conservation of energy, the right-hand side being the incoming and the left-hand side the outcoming and the absorbed energy.
>From the second estimate we deduce the following result.
\begin{proposition}\label{energie-demi-theo}
Let $(i\ep D - 2\overrightarrow{k}\cdot\overrightarrow{n})u^{in} \in L^2 (\R)$. 
If  $u \in H^1(\R_+\times \R)$ is a solution  to the problem (\ref{bbbb}) (\ref{CE0}) in the half-plane, it is unique. Moreover, we have the following stability estimate, with a constant $C$ independent of $\nu,$  $\mu$:
$$ \iint\limits_{\mathcal{D}}2\nu |u |^{2} + \int\limits_{\Gamma_0}  |\overrightarrow{k}\cdot\overrightarrow{n}| |u|^2  \leq C \int \limits_{\Gamma_0} |(i\ep D - 2\overrightarrow{k}\cdot\overrightarrow{n})u^{in}|^2.$$
\end{proposition}
{\bf Remark:} the stability result cannot be used to let the absorption tend to zero, since the norm on the left-hand side depends on the coefficient $\nu$.

In the sequel of this article, we assume that
\begin{equation}
\mu=0, \quad
\hbox{ and } \nu \hbox{ is a strictly positive constant} 
\end{equation}
 and we write 
$$2k_x g=i\ep k_y (k_x\Dy -k_y \Dx)u^{in}+2k_x u^{in}.$$

Problem  (\ref{bbb})(\ref{Cbord}) may be written under the form:
\begin{equation}\label{S-demi-2}
i(k_x \Dx + k_y \Dy)u + {\ep \over 2}(k_x^2 \Dyy -2k_x k_y \Dxy + k_y^2 \Dxx)u+i\nu u=0 \,on\,\R_+\times\R,
\end{equation} \nopagebreak
\begin{equation}\label{CE-demi-2}
i\ep k_y (k_x\Dy -k_y \Dx)u_{|x=0} +2k_x u_{|x=0} =g.
\end{equation}
We denote the Fourier variables of the variables $x$ and $y$ respectively by $\xi$ and $\eta$  and the Fourier transform in $x$ and $y$ by ${\cal F}_x$ and ${\cal F}_y$.
Using Fourier transformation, the following theorem gives an explicit solution to this problem. Denote
$$
 R_- (i\eta)=i{
 k_x  \eta \over k_y}-i{ k_x \over {\ep  k_y^2}}( 1 - \sqrt{1-2{\ep k_y  \eta \over k_x^2} + 2i \nu  {\ep k_y^2 \over k_x^2}}).
$$

Here and in the sequel, $\sqrt{~~}$ denotes the principal determination of the square root (its real part is positive); it is crucial that $\nu$ is strictly positive in order to define it without ambiguity.
\begin{theoreme} \label{demi-theo-2}
Let ${\cal S}'(\R)$ be the space of tempered distributions and $g \in {\cal S}'$. There exists a unique tempered distribution $u \in {\cal C}^\infty_x(\R_+, {\cal S}'_y(\R))$ to Problem (\ref{S-demi-2})(\ref{CE-demi-2}), which is:
\begin{equation}\label{ug}
{\cal F}_y (u;x,\eta)= {2{\cal F}_y (g;\eta) \over 1+ \sqrt{1-2{\ep k_y \eta \over k_x^2} + 2i\nu {\ep k_y^2 \over k_x^2}}}e^{R_- (i\eta) x}.
\end{equation}
It  satisfies also:
$$\biggl(\Dx - R_- (i\eta)\biggr) {\cal F}_y (u;x,\eta)=0.$$
                                \end{theoreme}
>From this theorem, we can also deduce the following corollaries.
\begin{coroll}\label{coroll-1}
If $g \in H^{-{1\over 2}}(\R)$ then the solution $u$ to Problem (\ref{S-demi-2})(\ref{CE-demi-2}) is in $ {\cal C}^{b}_x(\R_+, L^2_y (\R))$, and the following stability inequality stands for some constant $C$ not depending on the absorption factor $\nu$:
$$||u||_{L^{\infty}_x(\R_+, L^2_y (\R))} \leq C ||g||_{ H^{-{1\over 2}}(\R)}.$$
On a general manner, if $g \in H^s (\R)$, $s\in \R$, then $u \in {\cal C}^b_x (\R_+, H^{s+{1\over 2}}_y (\R))$ and we have the following inequality:
$$||u||_{L^{\infty}_x(\R_+, H^{s+{1\over 2}}_y (\R))} \leq C ||g||_{ H^{s}(\R)}.$$
\end{coroll}

The space chosen for the solution is the most convenient because the stability obtained does not depend on the absorption coefficient $\nu$. However, we can find existence and uniqueness of a solution in other spaces, provided other conditions on the initial data $g$. For instance:
\begin{coroll}\label{coroll-2}

\begin{enumerate}
\item

If ${{\cal F}_y (g;\eta) \over (1+ |\eta|^2)^{1\over 4}\sqrt{|{\cal R}e(R_- (i\eta))|} }\in L^2_\eta(\R)$ then $u \in L^2 (\R_+ \times \R).$ 
\item
If ${{\cal F}_y (g;\eta) (1+ |\eta|^2)^{s\over 2} \over \sqrt{|{\cal R}e(R_- (i\eta))|}} \in L^2_\eta(\R)$ with $s>0$ then $u \in L^\infty_y (\R, L^2_x (\R_+))$. In a general manner, if ${{\cal F}_y (g;\eta) (1+ |\eta|^2)^{s\over 2} |R_- (i\eta)|^m \over \sqrt{|{\cal R}e(R_- (i\eta))|}} \in L^2_\eta(\R)$ then $u \in L^\infty_y (\R, H^m_x (\R_+)).$
\end{enumerate}
\end{coroll}
Notice that when $\nu\to 0$, ${\cal R}e(R_- (i\eta)_{\nu=0})$ values zero in a half-line, then with these assumptions $u$ is not uniformly bounded according to the absorption coefficient $\nu$.

\ 

\begin{proof}
The proof is based on Fourier transforms of the problem, taken successively in $y$ and in $x.$ The factorization of the equation is the essential step of the proof, and will also provide us tools to define a proper transparent boundary condition for the quadrant problem.

Let $u$ be a solution of Problem (\ref{S-demi-2})(\ref{CE-demi-2}) and $v$ the extension of $u$ by zero in the whole space: $v(x,y)= u(x,y) \3 $. By introducing formally the function  $v$ in Equation (\ref{S-demi-2}) we get:
$$i\overrightarrow{k}\cdot\overrightarrow{\mathbf{\nabla }}v+\frac{\epsilon }{2}%
\Delta _{\bot }v+i\nu v= \biggl(\bigl(ik_x -{ \ep k_y \over 2}(2k_x \Dy -  k_y \Dx)\bigr)u(0,y)\biggl){\delta}_{x=0} + {\ep k_y^2 \over 2}u(0,y){\delta}^{'}_{x=0}.
$$
The meaning of the term $\Dx u(0,y)$ is given by the entrance boundary condition (\ref{CE-demi-2}) :
$$i\overrightarrow{k}\cdot\overrightarrow{\mathbf{\nabla }}v+\frac{\epsilon }{2}%
\Delta _{\bot }v+i\nu v=  ik_x g(y)\delta_{x=0} -{\ep k_y \over 2}\bigl( k_x \Dy u(0,y){\delta}_{x=0} - k_y u(0,y){\delta}^{'}_{x=0}\bigr).$$
Assuming that $u\in {\cal C}(\R_+, {\cal S}'(\R))$, we are allowed to take the Fourier transform of this expression. Let us define $P_{\nu}(X,Y)$ as the polynomial which characterizes the differential operator of the equation, \emph{i.e}:
 $$P_\nu(\Dx, \Dy)= i(k_x \dx + k_y \dy) + \frac{\epsilon }{2}(k_y^2 \Dxx-2k_x k_y \Dxy +k_x^2 \Dyy)+i\nu.  $$
If we write $u_0(y)=u(0,y)$, the Fourier transform in $y$ of the equation in $v$ reads:
$$P_\nu(\Dx, i\eta){\cal F}_y (v;x,\eta) =
 {\ep k_y^2 \over 2}\biggl\{\biggl({2i k_x \over \ep k_y^2}{\cal F}_y(g;\eta) - i {k_x\over k_y} \eta {\cal F}_y( u_0;\eta)\biggr)\delta_{x=0} + {\cal F}_y(u_0;\eta)\delta^{'}_{x=0}\biggr\}.$$
We can write 
\begin{equation}
P_\nu(\Dx,i\eta)={\ep k_y^2 \over 2}\biggl(\Dx -R_+ (i\eta)\biggr)\biggl(\Dx - R_- (i\eta)\biggr)
\label{fac1}
\end{equation}
where we define:
$$
R_\pm (i\eta)=i{
 k_x   \over k_y}\eta -i{ k_x \over {\ep  k_y^2}}\biggl( 1 \pm \sqrt{1-2{\ep k_y  \eta \over k_x^2} + 2i \nu  {\ep k_y^2 \over k_x^2}}\biggr).
$$
 Thus:
$$
\biggl(\Dx -R_+(i\eta)\biggr)\biggl(\Dx - R_-(i\eta)\biggr){\cal F}_y (v;x,\eta)=
$$\begin{equation}\label{new}
 \biggl({2i k_x \over \ep k_y^2}{\cal F}_y(g;\eta) - i {k_x\over k_y} \eta {\cal F}_y( u_0;\eta)\biggr)\delta_{x=0} + {\cal F}_y(u_0;\eta)\delta^{'}_{x=0}.
\end{equation}

We have reduced the problem in finding a unique acceptable solution for this ordinary differential equation (\ref{new}).
Its Fourier transform in $x$ reads:
$$\biggl(i\xi -R_+(i\eta)\biggr)\biggl(i\xi - R_-(i\eta)\biggr){\cal F}_x {\cal F}_y (v;\xi,\eta)= {2i k_x \over \ep k_y^2}{\cal F}_y(g;\eta) - i ({k_x\over k_y} \eta -\xi) {\cal F}_y( u_0;\eta).$$
Since ${\cal R}e \bigl(i\xi - R_\pm (i\eta)\bigr) \neq 0$, we can divide each side of Equation (\ref{new}) by ${2\over \ep k_y^2}P_\nu$ (because $1 \over P_\nu$ is an infinitely derivable function with at most polynomial growth at infinity, we can multiply it with a tempered distribution) and write:
$${\cal F}_x {\cal F}_y (v;\xi,\eta) = {\alpha_{1/2}^+ (\eta)\over i\xi- R_+(i\eta)} + {\alpha_{1/2}^-(\eta) \over i\xi - R_-(i\eta)},$$
where:
$$\alpha_{1/2}^+(\eta)= {R_+(i\eta) -i{k_x\over k_y}\eta \over R_+(i\eta) - R_-(i\eta)}{\cal F}_y(u_0;\eta) + {2ik_x  \over \ep k_y^2}{{\cal F}_y(g;\eta)\over R_+(i\eta) - R_-(i\eta)},$$
$$\alpha_{1/2}^-(\eta)=- {R_-(i\eta) -i{k_x\over k_y}\eta \over R_+(i\eta) - R_-(i\eta)}{\cal F}_y(u_0;\eta) - {2i k_x \over \ep k_y^2} {{\cal F}_y(g;\eta)\over R_+(i\eta) - R_-(i\eta)}.$$

If $\theta \in \C \backslash \R$, we can verify immediately that:
\begin{equation}\label{l1}
{1\over i\xi-\theta} = \left\{\begin{array}{ll}
 {\cal F}_x (\3 e^{\theta x};\xi) & \mbox{ if ${\cal R}e(\theta )<0$}\\
 -{\cal F}_x (\5 e^{\theta x};\xi) & \mbox{ if ${\cal R}e(\theta )>0$}
 \end{array}\right.\end{equation}
Applied to $\theta = R_+ (i\eta)$ and $\theta = R_- (i\eta),$ since ${\cal R}e (R_+)=-{\cal R}e (R_-) >0,$ we get:
 $${\cal F}_y(v;x,\eta) = -\alpha_{1/2}^+ (\eta) e^{R_+(i\eta)x}\5 + \alpha_{1/2}^-(\eta)e^{R_-(i\eta)x}\3,$$
which proves that  $v \in {\cal C}\bigl(\R_+, {\cal S}'(\R)\bigr)$.
Since $v(x<0,y)=0$, this equality implies $\alpha_{1/2}^+ (\eta)=0$. We take the limit for $x \stackrel{>}{\to}0$, and find ${\cal F}_y (u; 0, \eta)= \alpha_{1/2}^- (\eta).$ It gives Formula (\ref{ug}) and ends the proof of Theorem \ref{demi-theo-2}.

\end{proof}

Notice that with Formula (\ref{ug}), we easily calculate the value of the derivative $\overrightarrow{k}\cdot \nabla u$. As soon as $u$ is regular enough to have its Fourier transform in $y$ decreasing rapidly, we can develop it asymptotically in $\ep$ and $\nu$, and find: $\overrightarrow{k}\cdot \nabla u= O(\ep +\nu).$ It is a way to check \emph{a posteriori} the validity of our first approximation. 

\

{\bf Proof of the Corollaries.}

To show the stability result of Corollary \ref{coroll-1}, we use the following lemma:
                        \begin{lemma}
There exists a constant $C >0$, depending only on the geometry of the problem (\emph{i.e.} of $\ep$, $k_x$ and $k_y$) and not depending on $\nu$, such that:
$${1 + |\eta|^2 \over |1 +\sqrt{1 - {2\ep k_y\over k_x^2} \eta + 2i\ep \nu {k_y^2\over k_x^2}}|^4} \leq {C^4 \over 16}.$$
                        \end{lemma}
                \begin{proof}
Let us denote $X= 1 - {2\ep k_y\over k_x^2} \eta$ and $N=2\ep \nu {k_y^2\over k_x^2}$. Then:
$$B(\eta)=|1 + \sqrt{X +iN}|^2 =|1 + e^{i{\pi \over 4}}\sqrt{N-iX}|^2
$$
$$= 1 + \sqrt{X^2 + N^2} +2(X^2 + N^2)^{1\over 4} cos({\pi\over 4} - {Arctg{X\over N}\over 2}).$$
The cosine is nonnegative, so we get:
$$B(\eta)^2 \geq (1 + \sqrt{X^2 + N^2})^2 \geq 1 + X^2. 
$$
Since $1+|\eta|^2=1+a^2(1-X)^2$ with $a={k_x^2\over 2\ep k_y}$, we study the polynomial $P(X) = c(1+ X^2)-(1+a^2(1-X)^2)$, where $c={C^4\over 16}$. We look for a $c>0$ such that $\forall \; X \in \R$, $P(X)\geq 0$. 
$$P=X^2(c-a^2)+2a^2X-(1-c+a^2),$$
and we pick $c= 2a^2+1$ (for instance). We can take any $C \geq 2\sqrt{1 + {k_x^2 \over \ep |k_y|}}.$
\end{proof}

We now integrate in $\eta$ the square modulus of Equation (\ref{ug}). 

Since $|e^{R_-(i\eta)x}|= e^{{\cal R}e(R_-(i\eta))x} <1$, we get:
$$||u(x,.)||_{L^2(\R)}^2 \leq    C^2 \int\limits_{-\infty}^{\infty} {|{\cal F}_y(g;\eta)|^2 (1+|\eta|^2)^{-{1\over 2}}d\eta }= C^2 ||g||_{H^{-{1\over 2}}(\R)}^2. $$
We have proved the stability inequality, and the fact that the solution $u$ is in the space $L^\infty_x (\R_+, L^2 (\R))$. The constant $C$ is of order $1\over \sqrt{\ep}$. The proof is the same for $g\in H^s (\R).$

In order to prove 
Corollary \ref{coroll-2}, we first show that if $s>0$ and ${{\cal F}_y (g;\eta) (1+ |\eta|^2)^{s\over 2} \over \sqrt{|{\cal R}e(R_- (i\eta))|}} \in L^2_\eta(\R)$ then $u \in L^\infty_y (\R, L^2_x(\R_+))$. We have:
$$\int\limits_0^\infty |u(x,y)|^2 dx = \int\limits_0^\infty dx | \int {e^{i\eta y} d\eta \over 2\pi} {2{\cal F}_y(g;\eta)\over 1 + \sqrt{1 - {2 \ep k_y\over k_x^2}\eta + 2i\ep \nu {k_y^2\over k_x^2}}} e^{R_- (i\eta)x}|^2.$$
We can apply Cauchy-Schwartz inequality:
\begin{equation*}
\int\limits_0^\infty |u(x,y)|^2 dx \leq C^2 \int\limits_0^\infty dx \biggl(\int d\eta |{\cal F}_y (g;\eta)|^2 (1 + |\eta|^2)^{s} e^{2 {\cal R}e \bigl(R_- (i\eta)\bigr)x} \biggr) \biggl(  \int {d\eta (1+ |\eta|^2)^{-{1\over 2}-s}} \biggl).
\end{equation*}
it gives, for some constant $C>0$:
$$\int\limits_0^\infty |u(x,y)|^2 dx \leq A  \int  
d\eta  { |{\cal F}_y (g;\eta)|^2 (1 + \sqrt{|\eta|})^{2s}  \over 2 |{\cal R}e \bigl(R_- (i\eta)\bigr)| }.$$
This shows that $u \in L^\infty_y (\R, L^2_x(\R_+))$ if we suppose for instance ${{\cal F}_y (g;\eta) (1 + \sqrt{|\eta|})^s \over \sqrt{|{\cal R}e(R_- (i\eta))|}} \in L^2_\eta(\R).$ To generalize at $H^m_x (\R_+)$, we replace ${\cal F}_y (g)$ by $R_- (i\eta)^m {\cal F}_y (g;\eta)$.

                        \section{Case of the Quadrant}
\label{quart}

We still consider Problem (\ref{bbb})(\ref{Cbord}) but restricted to the quadrant $\{x,y\geq 0\}$. To find a \emph{transparent} or an \emph{absorbing} boundary condition (cf. Section \ref{intro-2}) on the boundary $\{y=0\}$, 
we formally factorize the differential operator of the advection-Schr\"odinger equation as follows: 
\begin{equation}  
P_\nu = i \overrightarrow{k}\cdot \overrightarrow{\mathbf{\nabla }}+ \frac{\epsilon }{2} \Delta _{\bot }u+i\nu=0= \ep {k_x^2 \over 2}\bigl(\Dy -A_+(\Dx)\bigr)\bigl(\Dy -A_-(\Dx)\bigr).
\label{fact}
\end{equation} 
where
$$A_- (\Dx)= {k_y \over k_x} \Dx -i {k_y\over \ep k_x^2}\bigl(1 - \sqrt{1+{2i\ep k_x\over k_y^2} \Dx +2i\ep\nu {k_x^2\over k_y^2}}\bigr).$$
It means that $A_+(X),A_-(X)$ are the roots of the polynomial $P_{\nu}(X,Y)$ considered as a polynomial of the variable $X$ (see Appendix \ref{racine-2} for a recall on the definition of a fractional derivative, or \cite{Schwartz}, chapter VI.5 for a complete construction).

The aim of this section is to show that the following condition:
$$\bigl(\Dy - A_+ (\Dx)\bigr)u_{|y=0}=0 \;\;\;\forall \;\; x>0 $$
is a \emph{transparent} boundary condition if $k_y >0$ and an \emph{absorbing} boundary condition if $k_y<0$. The proofs will make appear more clearly why we were lead to such a choice for a boundary condition.

\bigskip

Let us consider the following problem:
\begin{equation}\label{S-quart}
i(k_x \Dx + k_y \Dy)U + {\ep \over 2}(k_x^2 \Dyy -2k_x k_y \Dxy + k_y^2 \Dxx)U+i\nu U=0 \;\; \forall\;x>0,\;\; y>0,
\end{equation}
\nopagebreak
\begin{equation}\label{CE-quart}
i\ep k_y (k_x\Dy -k_y \Dx)U_{|x=0} +2k_x U_{|x=0} =2k_x g_+\;\; \forall \;\; y>0,
\end{equation}
\nopagebreak
\begin{equation}\label{CT}
\Dy U_{|y=0} - A_+ (\Dx) (U_{|y=0}) = 0 \;\; \forall x>0,
\end{equation}
where we take $g_+ \in H^{-{1\over 2}}(\R_+)$, $Supp(g_+)\subset \R_+^*$.
 We can rewrite Condition (\ref{CT}) under the following form:
\begin{equation*} \label{CTbis}
i\ep k_x (k_y\Dx -k_x \Dy)U_{|y=0} + k_y (1+\sqrt{1+{2i\ep k_x\over k_y^2}\Dx+2i\ep\nu{k_x^2\over k_y^2}})U_{|y=0}=0 \;\;\; \forall \;\; x>0.
\end{equation*}
The first step is to give an extended meaning of the boundary condition (\ref{CT}), for functions whose derivatives at the boundary $\{y=0\}$ may not be defined. The second step is
to show that Condition (\ref{CT}), considered in this extended meaning, is either a transparent or an absorbing boundary condition for our problem, according to the sign of $k_y$.  

Let us denote $u$  the solution of the half-space problem (\ref{S-demi-2})(\ref{CE-demi-2}) with an entrance boundary condition equal to $g=g_+\4$: $g \in H^{-{1\over 2}}(\R).$ It is given by Formula (\ref{ug}):
$${\cal F}_y (u;x,\eta)={2{\cal F}_y(g;\eta)\over 1 + \sqrt{1 - {2 \ep k_y\over k_x^2}\eta + 2i\ep\nu {k_y^2\over k_x^2}}} e^{R_- (i\eta)x}.$$
We show the two following theorems.
\begin{theoreme}\label{CT-extension}
Let $g_+$ be a given function such that $Supp(g_+)\subset \R_+^*$ and $g_+ \in H^{-{1\over 2}} (\R_+).$ 
Let $U\in  {\cal C}^b_x(\R_+, H^{{3\over 2}+s}_y (\R_+))$, with $s>0$, a function satisfying Equations (\ref{S-quart})(\ref{CE-quart}) with $\Dy U (x,0) \in H^{-{1\over 2}}_x (\R_+)$ and $U(x,0)\in H^{1\over 2}_x (\R_+)$. 
Then $U$ satisfies (\ref{CT}) \emph{iff} $U$ satisfies:
\begin{equation}\label{CT2}
U(x,y)\4={\cal F}_y^{-1} \biggl(\bigl\{{\h{K} (\eta) {\cal F}_y (U\4;0,\eta ) + \h{G}(\eta)\bigr\} e^{R_-(i\eta)x}} \biggr)\4,
\end{equation}
with:
$$\h{K} (\eta)=- {R_-(i\eta) -i{k_x\over k_y}\eta \over R_+(i\eta) - R_-(i\eta)} \;\;\;\; and \;\;\;\; \h{G}(\eta)= -{2ik_x \over \ep k_y^2} {{\cal F}_y(g;\eta)\over R_+(i\eta) - R_-(i\eta)}.$$

\end{theoreme}
With this result, we are able to extend the meaning of Equation (\ref{CT}) in the following way.
\begin{Def} 
We say that a function $U \in {\cal C}^b_x(\R_+, L^2_y (\R_+))$ is a solution of Problem  (\ref{S-quart})(\ref{CE-quart})(\ref{CT}) \emph{iff} it is a solution of Problem  (\ref{S-quart})(\ref{CE-quart})(\ref{CT2}). \label{def-CT}
\end{Def}
We are now able to state the main
                \begin{theoreme} \label{quart-theo1}
Let $g_+$ be a given function such that $Supp(g_+)\subset \R_+^*$ and $g_+ \in H^{-{1\over 2}} (\R_+).$

\noindent
i) Problem (\ref{S-quart})(\ref{CE-quart})(\ref{CT2}) admits a unique solution $U \in {\cal C}^b_x(\R_+, L^2_y (\R_+)) $.

\

\noindent
ii) Let $u\in {\cal C}^b_x(\R_+, L^{2}_y (\R))$ be the solution of the half-space problem (\ref{S-demi-2})(\ref{CE-demi-2}) with an incoming boundary condition $g=g_+ \4$.
We have an explicit formula for $U$ in terms of $u$, with $u_0(y) =u(0,y)$:
\begin{equation}\label{U=u}
U(x,y)\4={\cal F}_y^{-1} \biggl(\bigl\{\h{K}(\eta){\cal F}_y (u_0\4;\eta)+\h{G}(\eta)\bigr\}e^{R_-(i\eta)x}\biggr)\4,
\end{equation}
where $\h{K}=-{R_- -i{k_x\over k_y}\eta \over R_+ -R_-}$ and $\h{G}=-{2ik_x\over \ep k_y^2}{{\cal F}_y(g)\over R_+-R_-}$. 

\

\noindent
iii) If $k_y>0$, then $U=u_{|y\geq 0}$. 

\

\noindent
iv) If $k_y <0$, and if we take $g^A(y)=h(y-A)$ with $A>0$, denoting $u^A, U^A$ the corresponding solutions respectively in the half-space and in the quadrant, we get:
$$\lim\limits_{A\to +\infty} ||(u^A-U^A)\4||_{L^\infty(\R_+, L^2_y(\R))}=0.$$
                \end{theoreme}

                \subsection{Fourier transform of the problem.}

\label{Fourier-quart}

 Let $U$ be a solution of Problem (\ref{S-quart})(\ref{CE-quart})(\ref{CT}), $U\in {\cal C}^b_x (\R_+, H^{{3\over 2}+s}(\R_+))$ and $V$ the extension of $U$ by zero in the whole space: $V(x,y)= U(x,y) \3 \4 $. We calculate what the problem means for $V$, and find:
$$\begin{array}{l}
\displaystyle{
i \overrightarrow{k}\cdot\overrightarrow{\mathbf{\nabla }}V+\frac{\epsilon }{2}%
\Delta _{\bot }V+i\nu V=- \ep k_x k_y U(0,0)\delta_{x=0}\delta_{y=0} +}\\ \displaystyle{
\biggl(\bigl(ik_x -{ \ep k_y \over 2}(2 k_x \Dy -  k_y \Dx)\bigr)U(0,y)\biggr)\4 {\delta}_{x=0} + {\ep k_y^2 \over 2}U(0,y)\4{\delta}^{'}_{x=0} +} \\ \displaystyle{
\biggl(\bigl(ik_y -{ \ep k_x \over 2}(2k_y \Dx -  k_x \Dy)\bigr)U(x,0)\biggr)\3{\delta}_{y=0} + {\ep k_x^2 \over 2}U(x,0)\3{\delta}^{'}_{y=0},}
\end{array}
$$
which gives, with the two boundary conditions (\ref{CE-quart})(\ref{CT}):
\begin{eqnarray}\label{new2}
\begin{array}{l}
\displaystyle{
i\overrightarrow{k}\cdot\overrightarrow{\mathbf{\nabla }}V+\frac{\epsilon }{2}%
\Delta _{\bot }V+i\nu V=-\ep k_x k_y U(0,0)\delta_{x=0}\delta_{y=0}+  } \\
\displaystyle{
\biggl( i k_x  g \delta_{x=0} - {\ep k_y \over 2}\bigl( {k_x} \Dy U_{|x=0} \delta_{x=0} -k_y U_{|x=0}\delta^{'}_{x=0}\bigr)
\biggr)\4+  }\\
{
 \biggl(i{k_y\over 2}\bigl(1-\sqrt{1+{2i\ep k_x\over k_y^2}\Dx +2i\ep\nu{k_x^2\over k_y^2}}\bigr) U_{|y=0}\delta_{y=0} -{ \ep k_x \over 2}\bigl( k_y \Dx U_{|y=0} \delta_{y=0}
-k_x U_{|y=0}\delta^{'}_{y=0}\bigr)\biggr)\3 .}
\end{array}
\end{eqnarray}
To take the Fourier transform of this expression, we have to proceed carefully with the derivatives of $U$. Indeed, we can write:
$${\cal F}_x \bigl((\Dx U)\3;\xi,0\bigr)=i\xi {\cal F}_x \bigl(U\3; \xi,0\bigr)
-U(0,0),$$
and the equivalent formula for ${\cal F}_y \bigl((\Dy U)\4;0,\eta\bigr)$. On the contrary, one directly gets:
$${\cal F}_x \biggl(\bigl(\sqrt{1+{2i\ep k_x\over k_y^2}\Dx +2i\ep\nu{k_x^2\over k_y^2}} U_{|y=0}\bigr)\3;\xi,0\biggr)=\sqrt{1-{2\ep k_x\over k_y^2}\xi +2i\ep\nu{k_x^2\over k_y^2}}{\cal F}_x \bigl(U_{|y=0}\3\bigr).$$
We now take the Fourier transform of Formula (\ref{new2}). The terms depending on $U(0,0)$ cancel one another, and we furthermore notice that:
$$i{k_y\over 2}(1-\sqrt{1-{2\ep k_x\over k_y^2}\xi +2i\ep\nu{k_x^2\over k_y^2}})-{i\ep k_x k_y \over 2}\xi = -{\ep k_x^2\over 2} A_- (i\xi).$$
Finally, using the polynomial 
$P_{\nu}(X,Y)$, we get
\begin{equation*}
\begin{array}{l}
\displaystyle{
P_\nu(i\xi, i\eta){\cal F}_x {\cal F}_y (V;\xi,\eta) = } \\ \\
\displaystyle{ {\ep k_y^2 \over 2}\biggl({2i k_x \over k_y^2} {\cal F}_y(g;\eta) - i ({k_x\over k_y} \eta -\xi) {\cal F}_y( U\4;0,\eta)\biggr)  }
\displaystyle{+
{\ep k_x^2 \over 2} \bigl(i\eta -A_- (i\xi)\bigr){\cal F}_x(U\3;\xi,0).}
\end{array}
\end{equation*}
Dividing by $P_\nu$ written in one of the form (\ref{fact}) or (\ref{fac1}), the equation in $V$ reads:
\begin{equation}\label{general-quart}
{\cal F}_x {\cal F}_y (V;\xi,\eta) = {\alpha^+ (\eta)\over i\xi- R_+(i\eta)} + {\alpha^-(\eta) \over i\xi - R_-(i\eta)}+ {\beta^+ (\xi)\over i\eta- A_+(i\xi)} + {\beta^-(\xi) \over i\eta - A_-(i\xi)},
\end{equation}
where $\alpha^\pm$ and $\beta^\pm$ are given by
\begin{equation}\label{Fourier}
\begin{array}{l} 
\displaystyle{\alpha^+(\eta)= {R_+(i\eta) -i{k_x\over k_y}\eta \over R_+(i\eta) - R_-(i\eta)}{\cal F}_y(U\4;0,\eta) +{2ik_x \over \ep k_y^2} {{\cal F}_y(g;\eta)\over R_+(i\eta) - R_-(i\eta)},} \\ \\
\displaystyle{
\alpha^-(\eta)=- {R_-(i\eta) -i{k_x\over k_y}\eta \over R_+(i\eta) - R_-(i\eta)}{\cal F}_y(U\4;0,\eta) -{2ik_x \over \ep k_y^2} {{\cal F}_y(g;\eta)\over R_+(i\eta) - R_-(i\eta)},} \\ \\
\displaystyle{
\beta^+ (\xi)={\cal F}_x(U\3;\xi,0), }\\ \\
\displaystyle{
\beta^- (\xi)=0.}
\end{array}
\end{equation}

                         \subsection{Proof of Theorem \ref{CT-extension}.}

Let us take the inverse Fourier transform of Equation (\ref{general-quart}), as in Section \ref{demi-2} for the half-space problem:
$$\begin{array}{c}
U(x,y) \3\4= {\cal F}_y^{-1} \biggl(-\alpha^+ (\eta) e^{R_+(i\eta)x} \5 + {\alpha^-(\eta) e^{R_-(i\eta)x}\3} \biggr) \\
+ {\cal F}_x^{-1}\biggl(-\beta^+ (\xi) e^{A_+(i\xi)y}\8 + \beta^-(\xi)e^{A_-(i\xi)y}\4\biggr).
\end{array}
$$
If we multiply each side by $\3\4$ and use the fact that $\beta^- =0$, we obtain 
$$U(x,y) \3\4={\cal F}_y^{-1} \biggl({\alpha^-(\eta) e^{R_-(i\eta)x}} \biggr)\3 \4,$$
and denoting
$$\h{K} (\eta)=- {R_-(i\eta) -i{k_x\over k_y}\eta \over R_+(i\eta) - R_-(i\eta)} \;\;\;\; and \;\;\;\; \h{G}(\eta)= -{2ik_x \over \ep k_y^2} {{\cal F}_y(g;\eta)\over R_+(i\eta) - R_-(i\eta)},$$
we have $\alpha^- (\eta)=\h{K} (\eta) {\cal F}_y (U\4;0,\eta ) + \h{G}(\eta),$ and get Formula (\ref{CT2}).

It proves that  if $U \in {\cal C}^b_x (\R_+, H^{{3\over 2}+s}_y (\R_+))$ is a solution of Problem (\ref{S-quart})(\ref{CE-quart})(\ref{CT}), then it verifies Equation (\ref{CT2}).
Conversely, let  $U \in {\cal C}^b_x (\R_+, H^{{3\over 2}+s}_y (\R_+))$ be a solution of System (\ref{S-quart})(\ref{CE-quart})(\ref{CT2}), and $\Dy U(x,0)\in H^{-{1\over 2}}_x$ and $U(x,0)\in H^{1\over 2}_x$.  
Assuming no particular condition on the boundary $\{y=0\}$, we carry out the same computation
as in Section \ref{Fourier-quart}, and find:
$${\cal{F}}_{x,y}(U;\xi,\eta)={\alpha^+ (\eta)\over i\xi - R_+ (i\eta)} +{\alpha^- (\eta)\over i\xi - R_- (i\eta)} +{\beta^+ (\xi)\over i\eta - A_+ (i\xi)} +{\beta^- (\xi)\over i\eta - A_- (i\xi)},$$
with the same $\alpha^\pm$ than defined in (\ref{Fourier}) and with:
$$\beta^- (\xi)=-{ {\cal F}_x \biggl(\bigl(\Dy U-A_+ (\Dx)U\bigr)\3;\xi,0\biggr) 
\over A_+ (i\xi) - A_- (i\xi)}.$$
We take the inverse Fourier transform and multiply it by $\3\4$, to find:
$$U(x,y) \3\4 = {\cal F}_y^{-1} \biggl({\alpha^-(\eta) e^{R_-(i\eta)x}\3} \biggr)\4 + {\cal F}_x^{-1}\biggl(\beta^-(\xi)e^{A_-(i\xi)y}\4\biggr)\3. $$
Since $U$ verifies also Equation (\ref{CT2}), which writes:
$$U(x,y) \3\4={\cal F}_y^{-1} \biggl({\alpha^-(\eta) e^{R_-(i\eta)x}} \biggr)\3 \4,$$
we necessarily have:
$${\cal F}_x^{-1}\biggl(\beta^-(\xi)e^{A_-(i\xi)y}\4\biggr)\3=0 \;\;\; \forall x>0 \;\;, y>0. $$
The function $U$ being continuous in $y$, we can take the limit when $y$ tends to zero and obtain
${\cal F}_x^{-1}\bigl(\beta^-(\xi)\bigr)=0\;\;\; \forall x>0 $ (see \cite{Dym} and the theory of Hardy functions).

Since $\beta^-$ is the Fourier transform of a distribution null in $\R_-$ (see Appendix \ref{Hardy}), 
it implies that $\beta^-$ is the Fourier transform of a distribution whose support is included in $\{x=0\}$. Since $\beta^- \in L^2(\R)$, it implies $\beta^- =0$, which is equivalent to the Fourier transform in $x$ of the transparent boundary condition. 
We have proved Theorem \ref{CT-extension}.

Equation (\ref{CT}) can be defined only for functions in $H^{{3\over 2}+s}$ in the variable $y$. Theorem \ref{CT-extension} allows us now to extend it by Condition (\ref{CT2}), with Definition \ref{def-CT} given at the beginning of Section \ref{quart}.

                        \subsection{Proof of Theorem \ref{quart-theo1}.}

We recall Formula (\ref{CT2}):
$$U(x,y)\4={\cal F}_y^{-1} \biggl(\bigl\{{\h{K} (\eta) {\cal F}_y (U\4;0,\eta ) + \h{G}(\eta)\bigr\} e^{R_-(i\eta)x}} \biggr)\4.
$$
If we set $U_0 (y)= U(0,y)\4$, it gives:
\begin{equation}\label{alpha1}
U_0 (y) = {\cal F}_y^{-1} \biggl(\h{K} (\eta) {\cal F}_y (U_0;\eta ) + \h{G}(\eta)\biggr)\4.
\end{equation}
For these equations to make sense, it is sufficient to have $U\in  {\cal C}^b_{x} (\R_+, L^2_y (\R_+))$.
We prove the following lemma, from which we will easily deduce Theorem \ref{quart-theo1}.

                        \begin{lemma} \label{lemma-quart}

\noindent
i) Equation (\ref{CT2}) admits at most one solution $U\in {\cal C}^b_{x} (\R_+, L^2_y (\R_+))$; Equation (\ref{alpha1}) admits at most one solution $U_0 \in L^2_y (\R)$.

\noindent
ii) The solution $u\in {\cal C}^b_x (\R_+, L^2_y(\R))$ to Problem (\ref{S-demi-2})(\ref{CE-demi-2}) (with  $g=g_+ \4$) satisfies Equation (\ref{alpha1}). Moreover, $U=u_{|y>0}$ satisfies Equations (\ref{S-quart})(\ref{CE-quart}), and {\bf if $\bf{k_y >0}$ } it also satisfies Equation (\ref{CT2}).
                        \end{lemma}


\

\begin{proof}

i) Uniqueness of a solution of Equation (\ref{alpha1}) in $L^2(\R)$ implies that of a solution to Equation (\ref{CT2}) in ${\cal C}^b_x (\R_+, L^2_y (\R_+))$. Let us take $g=0$ and suppose that a function $U_0 \in L^2(\R)$ verifies:
$$U_0 (y) = {\cal F}_y^{-1} \biggl(\h{K} (\eta) {\cal F}_y (U_0;\eta)\biggr)\4.$$
Let $V_0 (y)={\cal F}_y^{-1} \biggl(\h{K}(\eta){\cal F}_y(U_0;\eta)\biggr) (y)$. 
The function $V_0$ verifies the following equation, because $U_0 = V_0 \4$:
$$V_0(y)={\cal F}_y^{-1} \biggl(\h{K}(\eta){\cal F}_y(V_0\4;\eta)\biggr)=\int\limits_0^\infty K(y-s)V_0(s)ds.$$
We separate $V_0$ into $V_+ = V_0\4$ and $V_- =V_0 \8$. Since $V_0 \in L^2 (\R)$, $V_\pm \in L^2 (\R_\pm)$ the functions $\h{V}_\pm := {\cal F}_y (V_\pm; \eta)$ belong to  the Hardy spaces $ {\cal H}^{2\pm}$ (see Appendix \ref{Hardy} for a recall on Hardy spaces).
The resulting equation writes:
$$\h{V}_+(\eta) (1-\h{K}(\eta))={1\over 2}\biggl(1+{1\over \sqrt{1-{2\ep k_y\over k_x^2}\eta + 2i\nu\ep{k_y^2\over k_x^2}}}\biggr)\h{V}_+ (\eta)=-\h{V}_- (\eta).$$
The idea is to find an ${\cal H}^{2-}$ function on the left-hand side, equal to a ${\cal H}^{2+}$ function on the right-hand side: since ${\cal H}^{2+}\cap {\cal H}^{2-}=0$, it will imply that both sides are null.
We use the fact that the function ${1\over \sqrt{1-{2\ep k_y\over k_x^2}\eta + 2i\nu\ep{k_y^2\over k_x^2}}}$ can be extended analytically by a uniformally bounded function, respectively: on $\R^{2-}=\{\eta -i\eta,\; \eta>0 \}$ if $k_y>0$, and on
$\R^{2+}$ if $k_y<0$. 
\begin{enumerate}
\item
$k_y >0:$ 
$1-\h{K}(\eta)$ is analytic and uniformly bounded in $\R^{2-}$, so
$(1-\h{K})\h{V}_+ \in {\cal H}^{2-}$, and is equal to $(-\h{V}_-) \in {\cal H}^{2+}$,
so $V_\pm=0$, so $V_0=0$. 

\item
$k_y <0.$
We write the previous equation in the following form:
$$\h{V}_+ (\eta)=-(1-\h{K}(\eta))^{-1}\h{V}_-(\eta).$$
$(1-\h{K}(\eta))^{-1}$ is analytic and uniformly bounded in $\R^{2+}$, so
$(1-\h{K})^{-1}\h{V}_- \in {\cal H}^{2+}$, and is equal to $(-\h{V}_+) \in {\cal H}^{2-}$, so $V_0=0$. 
\end{enumerate}

ii) Let $u\in {\cal C}^b_x (\R_+, L^2_y (\R))$ be the solution of Problem (\ref{S-demi-2})(\ref{CE-demi-2}) in the half-space, with an entrance boundary condition equal to $g\in H^{-{1\over 2}}(\R)$. 
We have seen in Section \ref{demi-2} that $u$ verifies:
$$u(x,y) \4={\cal F}_y^{-1} \biggl({\alpha^-_{1/2}(\eta) e^{R_-(i\eta)x}} \biggr)\4, $$
with $\alpha^-_{1/2}(\eta)= \h{K} (\eta)
{\cal F}_y(u;0,\eta) + \h{G}(\eta).$

This is almost Equation (\ref{CT2}): we have $\alpha^-=\alpha^-_{1/2} - \h{K} {\cal F}_y(u\8;0,\eta).$

 Let us distinguish the two cases, according to the sign of $k_y$.
\begin{enumerate}
\item $k_y >0.$
We have seen in Section \ref{demi-2} that $u$ verifies:
$${\cal F}_y(u;0,\eta) = {2{\cal F}_y(g\4;\eta)\over 1 + \sqrt{1 - {2 \ep k_y\over k_x^2}\eta + 2i\ep\nu {k_y^2\over k_x^2}}}.$$
Since the function ${1 \over 1 + \sqrt{1 - {2 \ep k_y\over k_x^2}\eta + 2i\ep\nu {k_y^2\over k_x^2}}}$ can be extended analytically by a uniformly bounded function in $\R^{2-}$, if $g\in L^2(\R_+)$ then ${\cal F}_y(u;0,\eta) \in {\cal H}^{2-}$, so $u(0,y)=u(0,y)\4$. 
This shows that in this case, $\alpha^-=\alpha^-_{1/2}$ so $u$ verifies Equation (\ref{CT2}).

\item $k_y <0.$
Since $u(0,.)\in L^2 (\R)$, $\h{K}(\eta){\cal F}_y(u\8;0,\eta) \in {\cal H}^{2+},$
 the function $u$ verifies the following equation:
$$u(0,y) \4={\cal F}_y^{-1} \bigl({\alpha^-(\eta) }\bigr)\4
+{\cal F}_y^{-1} \bigl(\h{K}(\eta){\cal F}_y(u\8;0,\eta)\bigr)\4.$$
Since  $\h{K}(\eta){\cal F}_y(u\8;0,\eta) \in {\cal H}^{2+},$ it is the Fourier transform of a function null in $\R_+$, so the last member of the equation values zero. This shows that $u(0,y)$ verifies Equation (\ref{alpha1}).
\end{enumerate}
\end{proof}

\

Let us deduce Assertions \emph{i)}, \emph{ii)} and \emph{iii)} of Theorem \ref{quart-theo1} from Lemma \ref{lemma-quart}. First, Lemma \ref{lemma-quart} proves the uniqueness of a solution to Problem (\ref{S-quart})(\ref{CE-quart})(\ref{CT}). In the case where $k_y >0$, we have nothing more to prove: the unique solution of Problem (\ref{S-quart})(\ref{CE-quart})(\ref{CT}) is the restriction to the quadrant of the solution of the half-space problem (\ref{S-demi-2})(\ref{CE-demi-2}).

\

Let us now take $k_y <0$. If $U$ exists, it verifies Equation (\ref{alpha1}) and we have seen that it is also the case of $u$, so $U (0,y)=u(0,y)_{|y\geq 0}$. Replacing $U(0,y)\4$ by $u(0,y)\4$ in Equation (\ref{CT2}) verified by $U$, this implies necessarily that:
$$U(x,y)\4={\cal F}_y^{-1} \biggl(e^{R_-(i\eta)x}\bigl\{\h{K}(\eta){\cal F}_y (u_0 \4;\eta)+\h{G}(\eta)\bigr\}\biggr)\4.$$
We now have to prove that the so-defined function $U$ is a solution of Problem (\ref{S-quart})(\ref{CE-quart})(\ref{CT}). First, it obviously verifies Equation (\ref{CT2}), so according to Definition \ref{def-CT} we have only to check Equations (\ref{S-quart}) and (\ref{CE-quart}). Let us define:
$$V:={\cal F}_y^{-1} \biggl(e^{R_-(i\eta)x}\bigl\{\h{K}(\eta){\cal F}_y (u_0\4;\eta)+\h{G}(\eta))\bigr\}\biggr).$$
We have $U=V\4$, so we can check Equations (\ref{S-quart}) (\ref{CE-quart}) on $V$ as well. We notice that $V=u-r$, with:
$${\cal F}_y (r;x,\eta)={\h{K}(\eta){\cal F}_y (u_0\8;\eta)}e^{R_- (i\eta)x}.$$ The function
$u$ verifies obviously Equations (\ref{S-quart})(\ref{CE-quart}), so we have only to check them on the remainder $r$.
Taking its global Fourier transform  and multiplying it by $P_\nu$, we find:
$$P_\nu (i\xi,i\eta){\cal F}_x {\cal F}_y (r\3;\xi,\eta)= (i\xi - R_+ (i\eta))\h{K}(\eta){\cal F}_y (u_0 \8;\eta)
.$$
Since ${\cal F}_y (u_0 \8;\eta) \in {\cal H}^{2+}$ and $k_y <0$, the right-hand side is the Fourier transform of a distribution null for $y>0$ (see Appendix \ref{Hardy} for a proof). Thus the remainder $r$ satisfies the Schr\"odinger equation (\ref{S-quart}) in the quadrant. Furthermore, at the boundary $x=0$, we have:
$${\cal F}_y (r;0,\eta)={\h{K}(\eta){\cal F}_y (u_0 \8;\eta)},$$
which is also  the Fourier transform of a distribution null for $y>0$. It proves that $r(0,y>0)=0$, so $r(0,.)_{|y>0}$ verifies the boundary condition (\ref{CE-quart}).

\

            {\bf Estimate of the difference between the solution in the quadrant and the solution in the half-space}
 
\

It only remains to prove Assertion \emph{iv)} of Theorem \ref{quart-theo1}. For the sake of notations simplicity, we skip here the indices and denote $u,\,U$ instead of $u^A, \,U^A.$
Let us assume that $k_y <0$ and denote by $u$ the solution of the half-space problem and by $U$ that in the quadrant. According to Equation (\ref{U=u}), we have:
$$U(x,y)={\cal F}_y^{-1} \biggl(e^{R_-(i\eta)x}[\h{K}(\eta){\cal F}_y (u_0\4;\eta)+\h{G}(\eta)] \biggr)\4,$$
and furthermore $u$ verifies the following equation:

$$u(x,y)\4={\cal F}_y^{-1} \biggl(e^{R_-(i\eta)x}[\h{K} (\eta) {\cal F}_y (u_0;\eta) + \h{G}(\eta)]\biggr)\4,$$
so by taking the difference, we obtain the relation :
\begin{equation}
(u-U)\4= {\cal F}_y^{-1}\biggl(e^{R_- (i\eta)x}\h{K}(\eta) {\cal F}_y (u\8;0,\eta)\biggr)\4.
\label{U-u}
\end{equation}

Let us now assume that $g(y)=h(y-A)$, $A>0$, $h\in H^{-{1\over 2}} (\R)$ and $Supp(h)\subset \R_+^*$. We can write:
$${\cal F}_y (u_0;\eta)={{\cal F}_y (h;\eta)e^{-i\eta A}\over 1 + \sqrt{1-2\ep{k_y\over k_x^2}\eta + 2i\nu\ep {k_y^2\over k_x^2}}}={\cal F}_y (H;\eta)e^{-i\eta A},$$
with a function $H \in L^2 (\R)$. We then have $u_0 = H(y-A)$, so $||u_0 \8||_{L^2}=||H_{|y\leq -A}||_{L^2},$ tends to $0$ when ${A\to +\infty}.$ We use Relation (\ref{U-u}) to estimate the difference between the half-space and the quarter plane solutions, since ${\cal R}e (R_- (i\eta))<0$:
$$||(u-U) (x,y)\4||_{L^\infty_x (\R_+, L^2_y(\R))} \leq ||\h{K}{\cal F}_y (u_0 \8)||_{L^2(\R)} \leq C ||H_{|y\leq -A}||_{L^2_y (\R)}.$$

\section*{Conclusion and perspectives}

A mathematical analysis has lead to an analytical form of the solution of the tilted paraxial equation in the case where the refraction index and the absorption coefficients are constant. We have proposed a convenient transparent/absorbing boundary condition for the problem on a quadrant, shown the well-posedness of the so-defined problem, and estimated the difference between the solutions on the half-space and on the quadrant.

However, this boundary condition is non-local, which implies much difficulty for its numerical treatment (see \cite{Arnold} for instance); up to now, the choice made in \cite{DDGS} to deal numerically with this boundary condition was to add an artificial absorption coefficient, as popularized in \cite{ber}. 

We have also noticed that Formula (\ref{ug}) allows us to check \emph{a posteriori} the validity of the formal asymptotic derivation of Equation (\ref{bbb}) from Helmholtz equation (\ref{base}). Complete estimates to justify it rigorously is a direction for future research (see also \cite{D} for a rigorous justification of the time-dependent problem). 

\

{\bf Acknowledgment.} The author expresses very grateful thanks to Fran\c{c}ois Golse and R\'emi Sentis for their precious help, guidance, ideas and corrections.

\appendix\section{Appendix}
\subsection{Fourier transforms of functions supported in $\R_+$}
\label{Hardy}

{\bf Background on Hardy classes}

For more details, we refer to \cite{Dym}.

\paragraph{Definition:}
let $h: \R^{2-}\to {\C}$. The function $h$ is said to be a \emph{Hardy function} if it verifies the two following properties:

i) h is analytic in $\R^{2-}= \{\w=a+ib, b<0\}$

ii) $\sup\limits_{b<0}\, ||h_b||^2_2 = \sup\limits_{b<0} \int\limits_{-\infty}^{\infty} |h(a+ib)|^2 da < \infty.$

\noindent
The so-defined space is called the \emph{Hardy space} and written as ${\cal H}^{2-}$.

\begin{theoreme} \label{H2-}
Let $h: \R^{2-}\to {\C}$. $h \in {\cal H}^{2-}$ \emph{iff} there exists a function $f\in L^2(\R_+)$ such that:
$$\forall \w \in \R^{2-}, \;\;\;\; h(\w) = \int\limits_0^\infty f(x) e^{i\w x} dx.$$
\end{theoreme}

{\bf Fourier transforms of functions that are identically zero in $\R_-$}

According to Theorem \ref{H2-}, the following kinds of tempered distribution $\h{f}$ are  Fourier transforms of tempered distributions $f$  supported in $\R_+$:

\begin{enumerate}
\item
$\h{f} \in {\cal H}^{2-} (\R)$: see Theorem \ref{H2-}.
\item
There exists a polynomial $P(i\xi)$ such that $\h{f}=P(i\xi)\h{g}(\xi)$, with $\h{g}$ the Fourier transform of a tempered distribution $g$, $Supp(g)\subset \R_+$: indeed, we then have $\h{f}={\cal F}_x (P(\Dx)g)$.

\item There exists a Hardy function $\phi \in {\cal H}^{2-}$ and a holomorphic function $\Psi,$ bounded on an open neighborhood of $\overline{\R^{2-}},$ such that $f=\phi\Psi.$

{\bf Examples:} 

\begin{enumerate}
\item the function $\Psi$ is a rational function with no pole in $\overline{\R^{2-}}\cup\{\infty\},$ \emph{i.e.} 
$$\Psi (z)= \sum\limits_{i=1}^N {\lambda_i \over z-\alpha_i},$$
with ${\cal I}m (\alpha_i) >0.$   

\item There exists $\alpha>0,$ $\beta\in\R,$ $s>0$ such that
$$\Psi (z)={1\over (1+i(\alpha z+\beta))^s}.$$ 
It is sufficient to prove this for $\Psi(z)={1\over (1+iz)^s}.$ On the one hand, we have:
$$1+iz \in \R_-^m\quad \equiv \quad iz \in (-\infty,-1)$$
$$1+iz \in \R_-^m\quad \equiv \quad z \in i(1,+\infty).$$ 
This implies that $\Psi$ is holomorphic on $\R^{2-}$ ($z^s = e^{s Log z}$ where $Log$ means the principal determination of the Logarithm).

On the other hand, $\Psi$ is bounded on a  neighborhood of $\R^{2-}:$ if ${\cal I}m(z)\leq 0$ we have
$$|\Psi(z)|={1\over |1+iz|^s}={1\over ((1-{\cal I}m(z))^2+{\cal R}e(z)^2)^{s/2}}\leq 1.$$ 

\end{enumerate} 

\item
The two previous cases imply that if $\h{g}\in {\cal H}^{2-}$ and $\h{f}=F(i\xi)\h{g}$, with $F$ a product of functions of the previous types, then $\h{f}$ is the Fourier transform of a distribution whose support is in $\R_+$.

\end{enumerate}

        \subsection{Definition of the square root of a differential operator}
\label{racine-2}

We recall here briefly how fractional derivatives are built, and refer to \cite{Schwartz}, chapter VI.5. for more details. For $a>0,$ we define
$$
Y_a(x)=\frac{(x_+)^a}{\Gamma(1+a)},
$$
where we have denoted $x_+= x \3.$ It is an homogeneous function of degree $a,$ which generalizes $\frac{x^n}{n!}.$ As for the polynomial case, one has
$$\forall \;a>1,\quad \frac{d Y_a}{dx}=Y_{a-1}.
$$
Taking this property as a definition, it allows us to define $Y_a$ by analytic continuation for $a\in \C.$ For the case $a\in -{\mathds N}^*,$ we notice that the fact that the pseudo-function  $Pf(x_+^a)$ is not defined is compensated by the poles of the $\Gamma$ function on $-{\mathds N}$ (see also \cite{Schwartz}, chapter II.3.) In particular, one has 
$$Y_0 (x)= sgn_+ (x),\quad Y_{-1} = \delta_{0},\quad Y_{-l}= \delta^{(l-1)}, \quad l\in{\mathds N}.$$
Since derivation corresponds to the convolution by $\delta'_0=Y_{-2},$ we can define $(\frac{d}{dx})^a,$ for functions supported in $\R_+,$ by the convolution by $Y_{-1-a}.$ It interpolates between $\frac{d}{dx}^0,$ which is the convolution by $\delta_0=Y_{-1},$ and $\frac{d}{dx},$ convolution by $\delta'_0=Y_{-2}.$ 
For the case $a=\frac{1}{2},$ we thus have:
$$
\left(\frac{d}{dx}\right)^{1/2}f=Y_{-3/2}\star
f=\frac{d}{dx}(Y_{-1/2}\star f).
$$
We can also have a formulation by the use of Fourier transforms:
$${\cal F} (Y_{-1/2})(\xi)=\frac{e^{-i.sign(\xi)\frac{\pi}{4}}}{\sqrt{|\xi|}}.$$
We obtain this formula by $\Gamma(\frac{1}{2})=\sqrt{\pi}$ and by passing to the limit when $\epsilon\to 0$ in
$$
\int_0^{+\infty}\frac{e^{-i\xi x}}{\sqrt{x}}dx=2\int_0^{+\infty}e^{-\epsilon y^2-i\xi y^2} dy
=\frac{\sqrt{\pi}}{\sqrt{\epsilon+i\xi}}.
$$
Finally:
$$
{\cal F}(\frac{d}{dx}(Y_{-1/2}\star f))=\frac{i\xi}{\sqrt{i\xi}} {\cal F}(f)(\xi).
$$
These two formulations give the two following equivalent ways to define properly the square root of the differential operator, which is, in our case, $-k_y^2-2i\ep k_x \Dx-2i\ep \nu k_x^2.$
\begin{enumerate}
\item
Using the Fourier transformation, we have:
$${\cal F}_x \biggl((\sqrt{-k_y^2-2i\ep k_x \Dx -2i\ep\nu k_x^2}u)\3; \xi\biggr)=$$
$$
 e^{-i{\pi \over 4}}\sqrt{-ik_y^2+2i\ep k_x \xi+2\ep\nu k_x^2} {\cal F}_x (u \3;\xi).$$

\item
Using  the expression $
\left(\frac{d}{dx}\right)^{1/2}f=Y_{-3/2}\star
f=\frac{d}{dx}(Y_{-1/2}\star f),$ we get:
$$\sqrt{-k_y^2-2i\ep k_x \Dx -2i\ep\nu k_x^2}(u)=
{\sqrt{2\ep k_x}e^{-i{\pi\over 4}+(i{k_y^2 \over 2\ep k_x} -\nu k_x) x}\over \sqrt{\pi}}\; {\Dx} \int\limits_0^x {u(s)e^{(-i{k_y^2 \over 2\ep k_x}+\nu k_x)s}\over \sqrt{x-s}}ds.$$

\end{enumerate}



%

\end{document}